# Robust estimation in finite population sampling

## Malay Ghosh[*1]

*University of Florida*

**Abstract:** The paper proposes some robust estimators of the finite population mean. Such estimators are particularly suitable in the presence of some outlying observations. Included as special cases of our general result are robust versions of the ratio estimator and the Horvitz-Thompson estimator. The robust estimators are derived on the basis of certain predictive influence functions.

## 1. Introduction

It is indeed a pleasure and privilege to contribute to this Festschrift honoring Professor P. K. Sen, a man whom I have long cherished as my friend, guide and philosopher. Professor Sen, for nearly a period of five decades, has made many profound contributions to the discipline of statistics. His research has encompassed every single area of statistical inference—parametric, semiparametric and nonparametric, and the theory that he has developed has found applications in many diverse areas of science. Indeed, he is one of the rare individuals in our profession who cannot just be identified with one localized area of statistics. The versatility of his research transcends any single narrowly focused topic, and the whole is by far bigger than the sum of the parts.

One of the many areas of interest of Professor Sen is the robustness of statistical procedures. His multiple authored or coauthored articles on the topic are very well summarized and unified in his 1996 classic treatise with Jureckova. The book provides a very comprehensive account of the subject with a fully developed asymptotic theory.

In this note, I will consider the robustness issue from a Bayesian perspective in the context of finite population sampling. Although, written within a Bayesian framework, the proposed estimators can also be viewed entirely from a model-based perspective. We introduce the notion of "predictive influence functions" as introduced by Johnson and Geisser [15, 16, 17], and obtain robust alternatives to a general class of Bayesian model-based estimators of the finite population mean, which includes in particular robust alternatives to the popular ratio estimators as well as Horvitz-Thompson estimators.

Section 2 of this paper introduces the concept of "predictive influence functions" based on a general divergence measure as introduced for example in Amari [1] and Cressie and Read [5]. The general divergence measure includes as a special case both the Kullback-Leibler and Hellinger divergence measures. The concept of "predictive influence functions" was first used by Johnson and Geisser [15, 16, 17] in a Bayesian context based specifically on the Kullback-Leibler divergence measure. Based on these influence functions, we have developed in this section, a general class of robust Bayesian estimators of the finite population mean. As special cases, we have found robust Bayesian alternatives to the ratio estimator as well as the Horvitz-Thompson estimator. In Section 3, we have obtained the mean squared error of these robust Bayes estimators purely from a frequentist criterion. Some final remarks are made in Section 4.

---

*Supported in part by NSF Grants SES-9911485 and SES-0631426.
[1]Department of Statistics, University of Florida, Box 118545, Gainesville, FL 32611-8545, USA, e-mail: ghoshm@stat.ufl.edu

Influence functions have a long and rich history in the statistics literature. Their importance in robust estimation is well-emphasized in Hampel [12], Huber [14], Hampel et al. [13] and many related papers. The predominant idea is to detect influential observations in terms of their effects on parameters, most typically the regression parameters. We have taken instead the predictive point of view. In finite population sampling, where the primary goal is to predict the "unseen" from the "seen", such a predictive approach seems quite natural. However, we have been able to point out a close connection between the proposed influence function, and the ones considered by Hampel, Huber and others. The main point is to control the effect of outlying observations for inference in finite population sampling. We have pointed out also the connection of the proposed robust ratio estimators to the corresponding estimators of Chambers [4] and Gwet and Rivest [10].

## 2. Development of estimators

Consider a finite population with units labeled $1, 2, \ldots, N$. A subset $s$ of $\{1, 2, \ldots, N\}$ is referred to as a sample. For simplicity, we consider only samples of fixed size $n$. Let $y_1, \ldots, y_N$ denote the characteristics of interest associated with the $N$ units in the population. Consider the hierarchical Bayesian model where conditional on $\theta$, $y_i \overset{\text{ind}}{\sim} N(\theta a_i, \sigma_i^2)$, $1, 2, \ldots, N$, where the $a_i$ and $\sigma_i^2$ are known constants, while the unknown $\theta$ has a uniform distribution on the real line. Without loss of generality, let $s = \{1, 2, \ldots, n\}$. Also, let $\bar{y}_w = \sum_{i=1}^n a_i \sigma_i^{-2} y_i / \sum_{i=1}^n a_i^2 \sigma_i^{-2}$, $\boldsymbol{a}_s = (a_1, \ldots, a_n)^T$, $\boldsymbol{a}_u = (a_{n+1}, a_{n+2}, \ldots, a_N)$, $\boldsymbol{y}_s = (y_1, \ldots, y_n)^T$, $\boldsymbol{y}_u = (y_{n+1}, y_{n+2}, \ldots, y_N)^T$ and $\boldsymbol{\Sigma}_u = \text{diag}(\sigma_{n+1}^2, \ldots, \sigma_N^2)$. It is shown in Ghosh and Sinha [9] that the posterior distribution of $\theta$ given $\boldsymbol{y}_s$ is $N(\bar{y}_w, 1/\sum_{i=1}^n a_i^2 \sigma_i^{-2})$, and the posterior predictive distribution of $\boldsymbol{y}_u$ given $\boldsymbol{y}_s$ is $N(\bar{y}_w \boldsymbol{a}_u, \boldsymbol{\Sigma}_u + \boldsymbol{a}_u \boldsymbol{a}_u^T / \sum_{i=1}^n a_i^2 \sigma_i^{-2})$. It is also shown that with the given model, the estimator of the finite population mean $\bar{y}_P = N^{-1} \sum_{i=1}^N y_i$ is given by

$$(1) \qquad \hat{\bar{y}}_P = N^{-1}[\sum_{i=1}^n y_i + \bar{y}_w \sum_{j=n+1}^N a_j].$$

In particular, if $a_i = x_i$, $\sigma_i^2 = \sigma^2 x_i$, $i = 1, \ldots, n$, then the resulting estimator of the finite population mean is the ratio estimator $(\bar{y}_s/\bar{x}_s)\bar{x}_P$, where $\bar{y}_s = n^{-1} \sum_{i=1}^n y_i$, $\bar{x}_s = n^{-1} \sum_{i=1}^n x_i$, and $\bar{x}_P = N^{-1} \sum_{i=1}^N x_i$. The choice $a_i = x_i$ and $\sigma_i^2 = x_i^2$ leads to the estimator $N^{-1}[\sum_{i=1}^n y_i + n^{-1} \sum_{i=1}^n (y_i/x_i)]$, an estimator introduced in Royall [19] and considered at length in Basu [2]. Finally, the choice $a_i = \pi_i$ and $\sigma_i^2 = \pi_i/(1 - \pi_i)$, where $\pi_i > 0$ for all $i = 1, \ldots, N$ and $\sum_{i=1}^N \pi_i = n$ leads to the Horvitz-Thompson estimator $N^{-1} \sum_{i=1}^n (y_i/\pi_i)$ for $\bar{y}_P$. It is instructive to view the estimator $\bar{y}_w$ of $\theta$ as a weighted average of the estimators $y_i/a_i$ of $\theta$, the weights being proportional to $a_i^2/\sigma_i^2$, the inverses of the variances of these estimators. However, the estimators $y_i/a_i$ are not necessarily outlier resistant. In particular, some of the components $y_i/a_i$ could be substantially different from the grand average $\bar{y}_w$. This may happen when for instance, a particular ratio $y_i/a_i$ is substantially smaller in magnitude than the other ratios $y_j/a_j$, but its variance $\sigma_i^2/a_i^2$ is much larger relative to the other variances $\sigma_j^2/a_j^2$, $(j \neq i)$. Since the estimators $y_i/a_i$ of $\theta$ are weighted inversely proportional to their variances in $\bar{y}_w$, the effect of this particular $y_i/a_i$ can be very insignificant compared to the other $y_j/a_j$ $(j \neq i)$, in finding $\bar{y}_w$. In order to control such outlying observations, as in Hampel [12], we first study the influence of $y_i/a_i$. To this end, we bring in the notion of "predictive influence functions" as introduced by Johnson and Geisser [15, 16, 17] based on the Kullback-Leibler (K-L) divergence measure. The influence function as considered here, however, is based on a general divergence measure introduced by Amari [1] and Cressie and Read [5]. This measure includes, but is not limited to the K-L or Bhattacharyya-Hellinger (B-H) (Bhattacharyya [3] and Hellinger [11]) divergence measure. For two densities $f_1$ and $f_2$, this general divergence measure is given by

$$(2) \qquad D_\lambda(f_1, f_2) = \frac{1}{\lambda(\lambda+1)} E_{f_1}\left[\left(\frac{f_1}{f_2}\right)^\lambda - 1\right].$$

The above divergence measure should be interpreted as its limiting value when $\lambda \to 0$ or $\lambda \to -1$. We may note that $D_\lambda(f_1, f_2)$ is not necessarily symmetric in $f_1$ and $f_2$, but the symmetry can always be achieved by considering $\frac{1}{2}[D_\lambda(f_1, f_2) + D_\lambda(f_2, f_1)]$. Also, it may be noted that as $\lambda \to 0$, $D_\lambda(f_1, f_2) \to E_{f_1}[\log \frac{f_1}{f_2}]$, while if $\lambda \to -1$, $D_\lambda(f_1, f_2) \to E_{f_2}[\log \frac{f_2}{f_1}]$. These are the two K-L divergence measures. Also $D_{-\frac{1}{2}}(f_1, f_2) = 4(1 - \int \sqrt{f_1 f_2}) = 2H^2(f_1, f_2)$, where $H(f_1, f_2) = \{2(1 - \int \sqrt{f_1 f_2})\}^{1/2}$, the B-H divergence measure. In the present context, we consider the divergence between the posterior predictive distribution of $\boldsymbol{y}_u$ given $\boldsymbol{y}_s$ and the posterior predictive distribution of $\boldsymbol{y}_u$ given $\boldsymbol{y}_s$ with one of the $y_i$, say, $y_k$, $k = 1, \ldots, n$ removed. To this end, we first state a general divergence result involving two multivariate normal distributions based on the general divergence measure as given in (2). The result is proved in Ghosh, Mergel and Datta (2006).

**Theorem 1.** *Let $f_1$ and $f_2$ denote the $N_p(\boldsymbol{\mu}_1, \boldsymbol{\Sigma}_1)$ and $N_p(\boldsymbol{\mu}_2, \boldsymbol{\Sigma}_2)$ pdf's respectively. Then*

$$\begin{aligned}
D_\lambda(f_1, f_2) &= \frac{1}{\lambda(\lambda+1)}[\exp\{\frac{\lambda(\lambda+1)}{2}(\boldsymbol{\mu}_1 - \boldsymbol{\mu}_2)^T \\
&\qquad\qquad \times ((1+\lambda)\boldsymbol{\Sigma}_2 - \lambda\boldsymbol{\Sigma}_1)^{-1}(\boldsymbol{\mu}_1 - \boldsymbol{\mu}_2)\} \\
&\quad \times |\boldsymbol{\Sigma}_1|^{-\frac{\lambda}{2}}|\boldsymbol{\Sigma}_2|^{-\frac{\lambda-1}{2}}|(1+\lambda)\boldsymbol{\Sigma}_2 - \lambda\boldsymbol{\Sigma}_1|^{\frac{1}{2}} - 1].
\end{aligned}$$

It follows from the above general result that the divergence between two normal distributions is a quadratic function of the difference of the two mean vectors. In the present context the difference in the mean vectors of the two posterior predictive distributions of $\boldsymbol{y}_u$ turns out to be a multiple of the square of $\bar{y}_w - \frac{\sum_{i=1}^n a_i\sigma_i^{-2}y_i - a_k\sigma_k^{-2}y_k}{\sum_{i=1}^n a_i^2\sigma_i^{-2} - a_k^2\sigma_k^{-2}}$, which on simplification reduces to a known multiple of the square of $y_k/a_k - \bar{y}_w$. Thus, one needs to control the residuals $y_k/a_k - \bar{y}_w$ for finding robust estimators of the finite population mean $\bar{y}_P$. However, in order to make these residuals scale-free, we consider the standardized residuals $(y_k/a_k - \bar{y}_w)/v_k$, where $v_k^2 = V(y_k/a_k - \bar{y}_w) = (a_k^2 \sigma_k^{-2})^{-1} - (\sum_{i=1}^n a_i^2 \sigma_i^{-2})^{-1}$.

It is instructive to find a connection between the proposed method of finding the influence functions, and those that are widely used in the robust statistics literature. Following the approach of Hampel [12], the influence function (IF) of the functional T at distribution function F is given by

$$\mathrm{IF}(x; T, F) = \lim_{t \downarrow 0} \frac{T((1-t)F + t\delta_x) - T(F)}{t},$$

for $x \in \mathcal{X}$ where this limit exists. Here $T(F)$ is the parameter of interest, and $\delta_x$ is the dirac-delta function. Thus when the parameter of interest is the population mean, namely $T(F) = \int x dF$, $X \sim F$, writing $F_t = (1-t)F + t\delta_{x_0}$,

$$T(F_t) = \int x d[(1-t)F + t\delta_{x_0}] = (1-t)\int x dF + tx_0.$$

for $x \in \mathcal{X}$ where this limit exists. Here $T(F)$ is the parameter of interest, and $\delta_x$ is the dirac-delta function. Thus when the parameter of interest is the population mean, namely $T(F) = \int x dF$, $X \sim F$, writing $F_t = (1-t)F + t\delta_{x_0}$,

$$T(F_t) = \int x d[(1-t)F + t\delta_{x_0}] = (1-t)\int x dF + tx_0.$$

Hence, $\mathrm{IF}(x_0; T, F) = x_0 - \theta$. In our case, conditional on $\theta$, $E(y_k/a_k) = \theta$, and the natural estimator of $\theta$ is $\bar{y}_w$. Hence, we estimate $y_k/a_k - \theta$ by $y_k/a_k - \bar{y}_w$. However, it is more appropriate to consider the scale-free residuals $r_k = (y_k/a_k - \bar{y}_w)/v_k$.

Based on these scale-free residuals, and writing $w_i = a_i^2 \sigma_i^{-2} / \sum_{i=1}^n a_i^2 \sigma_i^{-2}$, $i = 1, \ldots, n$, we propose the robust estimator

$$\hat{\theta}_R = \bar{y}_w + \sum_{i=1}^n w_i v_i [r_i I_{[|r_i| \le C]} + C I_{[r_i > C]} + (-C) I_{[r_i < -C]}]. \tag{3}$$

Consequently, the proposed robust estimator of the finite population mean $\bar{y}_P$ is given by

$$\bar{y}_P^{(R)} = N^{-1} [\sum_{i=1}^n y_i + \hat{\theta}_R \sum_{j=n+1}^N a_j]. \tag{4}$$

**Remark 1.** The proposed estimator of $\theta$ or of $\bar{y}_P$ is similar in spirit to the "limited translation estimator" of Efron and Morris [6, 7]. However, the present motivation of these estimators from the predictive influence function point of view is entirely new. We will address the question of choice of the constant $C$ in the next section.

**Remark 2.** In the special case when $a_i = \pi_i$ and $\sigma_i^2 = \pi_i^2/(1 - \pi_i)$, $\pi > 0$ and $\sum_{i=1}^N \pi_i = n$, $w_i = (1 - \pi_i)/\sum_{i=1}^n (1 - \pi_i))$ and $v_i = (1 - \pi_i)^{-1} - (\sum_{i=1}^n (1 - \pi_i))^{-1}$. The resulting estimator $\bar{y}_P^{(R)}$ of $\bar{y}_P$ is a robust alternative of the celebrated Horvitz-Thompson estimator.

**Remark 3.** Next in the case when $a_i = x_i$ and $\sigma_i^2 = \sigma^2 h(x_i)$, it follows that $w_i = [x_i^2/h(x_i)]/\sum_{i=1}^n [x_i^2/h(x_i)]$. In this case, our estimator is similar to the one of Chambers [4] except that Chambers used $\sigma_i$ rather than $v_i$ as the scaling factor. To see the difference between the two estimators in the special case of robust alternatives to the ratio estimator, that is, where $a_i = x_i$ and $\sigma_i^2 = \sigma^2 x_i$, the proposed scaling factor $v_i$ simplifies to $\sigma[x_i^{-1} - (\sum_{i=1}^{n} x_i)^{-1}]^{1/2}$. In contrast, Chamber's scaling factor is just $\sigma x_i^{-1/2}$. One interesting feature is that we bring in the notion of influence functions in the derivation of robust estimators in finite population sampling which could be potentially useful as a general approach to robust estimation in the model-based analysis of more complex surveys as well.

In the next section, we will find the frequentist mean squared error (that is conditional on $\theta$) of $\bar{y}_P^{(R)}$.

## 3. Mean squared error

We first find an expression for the mean squared error (MSE) of $\bar{y}_P^{(R)}$ as an estimator of $\bar{y}_P$. In the process, we have also made some observations about the choice of the constant $C$. To this we first prove the following theorem.

**Theorem 2.** *Under the assumed model, conditional on $\theta$,*

$$(5) \quad E(\bar{y}_P^{(R)} - \bar{y}_P)^2 = N^{-2}[\sum_{j=n+1}^{N} \sigma_j^2 + \{(\sum_{i=1}^{n} a_i^2 \sigma_i^{-2})^{-1} + (\sum_{i=1}^{n} w_i^2 v_i^2)g(C)\}(\sum_{j=n+1}^{N} a_j)^2],$$

*where $\Phi$ and $\phi$ denote respectively the N(0,1) df and pdf and*

$$(6) \quad g(C) = 2[(C^2+1)\Phi(-C) - 2C\phi(C)].$$

**Remark 4.** It may be noted that under the assumed model, the MSE of $\hat{\bar{Y}}_P$, the posterior mean of $\bar{Y}_P$, also the best unbiased predictor (best linear unbiased predictor without normality), is given by $N^{-2}[\sum_{j=n+1}^{N} \sigma_j^2 + (\sum_{i=1}^{n} a_i^2 \sigma_i^{-2})^{-1} \times (\sum_{j=n+1}^{N} a_j)^2]$. Thus, if the assumed model is true, the excess risk of the proposed robust estimator is given by $N^{-2}(\sum_{i=1}^{n} w_i^2 v_i^2)g(C)(\sum_{j=n+1}^{N} a_j)^2$. Noting that $g'(C) = 2[C\Phi(-C) - \phi(C)] < 0$ (Feller [8], page 166), it follows that $g(C)$ is decreasing in $C$. This is intuitively expected since larger the value of $C$ closer $\hat{\theta}_R$ is to $\hat{\theta}$. The constant $C$ will be chosen by setting an upper bound, say, $M$ to this excess risk, and then solving $C$ numerically by equating this excess risk to $M$. The choice of $M$ will be clearly left to the experimenter. The main idea is to seek a tradeoff between robustness against model failure and the maximum excess risk that one is willing to tolerate by proposing this robust estimator when the assumed model is true.

*Proof of Theorem 2.* Throughout, the calculations are done conditional on $\theta$. First from the independence of the $y_i$, $i = 1, \ldots, N$, for fixed $\theta$, it follows that

$$(7) \quad E(\bar{y}_P^{(R)} - \bar{y}_P)^2 = N^{-2}[\sum_{j=n+1}^{N} \sigma_j^2 + (\sum_{j=n+1}^{N} a_j)^2 E(\hat{\theta}_R - \theta)^2].$$

Next noting that $\sum_{i=1}^{n} w_i v_i r_i = \bar{y}_w - \bar{y}_w = 0$, from (3), one can alternately write $\hat{\theta}_R$ as

$$(8) \quad \hat{\theta}_R = \bar{y}_w - \sum_{i=1}^{n} w_i v_i[(r_i - C)I_{[r_i > C]} + (r_i + C)I_{[r_i < -C]}].$$

Next due to the independence of the $y_i$, $\text{Cov}(y_i/a_i - \bar{y}_w, \bar{y}_w) = (\sum_{i=1}^{n} a_i^2 \sigma_i^{-2})^{-1} - (\sum_{i=1}^{n} a_i^2 \sigma_i^{-2})^{-1}) = 0$. Hence, because of normality, $y_i/a_i - \bar{y}_w$, and accordingly $r_i$ is distributed independently of $\bar{y}_w$. Accordingly, it follows from (8) that

$$(9) \qquad E(\hat{\theta}_R - \theta)^2 = V(\bar{y}_w) + E[\sum_{i=1}^{n} w_i v_i \{(r_i - C)I_{[r_i > C]} + (r_i + C)I_{[r_i < -C]}\}]^2.$$

Since $r_i \sim \text{N}(0,1)$, it follows after some algebra that

$$(10) \qquad E[(r_i - C)I_{[r_i > C]} + (r_i + C)I_{[r_i < -C]}]^2 = 2[(C^2 + 1)\Phi(-C) - 2C\phi(C)].$$

Next, by the fact that for $i \neq k$, $(r_i, r_k) \overset{\mathrm{d}}{=} (r_i, -r_k) \overset{\mathrm{d}}{=} (-r_i, r_k) \overset{\mathrm{d}}{=} (-r_i, -r_k)$,(where $\overset{\mathrm{d}}{=}$ signifies "has the same distribution as"), one gets

$$E[\{(r_i - C)I_{[r_i > C]} + (r_i + C)I_{[r_i < -C]}\}\{(r_k - C)I_{[r_k > C]} + (r_k + C)I_{[r_k < -C]}\}]$$
$$= E[\{(r_i - C)(r_k - C) - (r_i - C)(r_k - C)$$
$$(11) \qquad -(r_i - C)(r_k - C) + (r_i - C)(r_k - C)\}I_{[r_i > C]I_{[r_k > C]}}] = 0.$$

The theorem follows now from (7) and (9)–(11). $\qquad \blacksquare$

## 4. Summary and conclusion

The paper proposes some robust estimators which can guard against outlying observations in connection with model-based inference in finite population sampling. In the process, new robust alternatives to the ratio estimators as well as the Horvitz-Thompson estimator are found. The mean squared errors of these model-based estimators are also obtained. Future work will encompass extension of these ideas to more complex surveys, for example in multistage stratified sampling, and also to address situations when there is wide departure from the assumed model.

**Acknowledgments.** Thanks are due to a referee for constructive comments.